\documentclass{amsart}

\usepackage{amsmath}
\usepackage{amscd}
\usepackage{amssymb}

\newcommand{\bC}{{\mathbb C}}
\newcommand{\bCP}{{\mathbb C}{\mathbb P}
}

\newcommand{\bZ}{{\mathbb Z}}
\newcommand{\C}{ {\mathbb C} }

\newcommand{\G}{G}
\newcommand{\Gn}{\G_n}
\newcommand{\Z}{ {\mathbb Z} }

\DeclareMathOperator{\diag}{diag}

\newtheorem{theorem}{Theorem}[section]
\newtheorem{theorem/definition}{Theorem/Definition}[section]

\newtheorem{proposition}{Proposition}[section]
\newtheorem{lemma}{Lemma}[section]

\newtheorem{corollary}{Corollary}[section]

\theoremstyle{remark}
\newtheorem{remark}{Remark}[section]

\theoremstyle{definition}
\newtheorem{example}{Example}[section]
\newtheorem{conjecture}{Conjecture}[section]

\begin{document}
\title
{Orbifold Hodge numbers of the wreath product orbifolds}
\author{Weiqiang Wang}
\address{Department of Mathematics\\
North Carolina State University\\ Raleigh, NC 27695}
\email{wqwang@math.ncsu.edu}
\author{Jian Zhou}
\address{Department of Mathematics\\
Texas A\&M University\\
College Station, TX 77843}
\email{zhou@math.tamu.edu}
\begin{abstract}
We prove that the wreath product orbifolds studied earlier by the
first author provide a large class of higher dimensional examples
of orbifolds whose orbifold Hodge numbers coincide with the
ordinary ones of suitable resolutions of singularities. We also
make explicit conjectures on elliptic genera for the wreath
product orbifolds.
\end{abstract}

\maketitle

\vspace{.3cm}

\noindent{Key words:} Wreath products, (orbifold) Hodge numbers,
Hilbert schemes, elliptic genera.

\vspace{.3cm} \noindent{2000 Mathematics Subject Classifications:}
14J, 58J.

\section{Introduction}
In the study of orbifold string theory, Dixon, Harvey, Vafa and
Witten \cite{DHVW} introduced the notion of orbifold Euler number
for a smooth manifold $Y$ acted on by a finite group $G$ and
raised the question on the existence of a resolution of
singularities of $M/G$ whose ordinary Euler number coincides with
the orbifold Euler number. The orbifold Euler number has
subsequently been interpreted as the Euler number for equivariant
$K$-theory, cf. Atiyah-Segal \cite{Ati-Seg}. The notion of
orbifold Euler number has been further refined to give rise to the
notion of orbifold Hodge numbers \cite{Vafa, Zas} and more
generally the stringy Hodge numbers \cite{Bat-Dai}. The orbifold
Hodge numbers of an orbifold are then conjectured to coincide with
the ordinary Hodge numbers of a suitable resolution of the
orbifold. For recent related development, see \cite{Bat-Dai, Bat,
DL, R} and the references therein.

A well-known series of examples with such a property is provided
by the symmetric product of a surface which admits a resolution of
singularities given by the Hilbert scheme of points. In this case
the orbifold Euler number calculated by Hirzebruch and H\"ofer
\cite{HH} matches with the Euler number of the Hilbert scheme
found earlier by G\"ottsche \cite{G1}.  It is further shown by
G\"ottsche \cite{G2} that the orbifold Hodge numbers matches with
the Hodge numbers of the Hilbert scheme calculated by G\"{o}ttsche
and Soergel \cite{Got-Soe} (also see Cheah \cite{Che}). The same
method has been used by the second author \cite{Zh}, also for the
calculation for higher dimensional complex
manifolds.\footnote{Seeing the math review (99c:14022) of \cite{G2} but not
the paper itself when he was writing \cite{Zh}, the second author
got the wrong impression that G\"{o}ttsche proved his result by
establishing the strong McKay correspondence for symmetric
products. This misunderstanding has been clarified when we
actually looked into the paper \cite{G2} during the preparation of
the present paper.}

The wreath product orbifolds, which are generalizations of the
symmetric products, were shown by the first author \cite{W} (also
see \cite{FJW, W2}) to have deep connections with Hilbert schemes
of surfaces and vertex representations of infinite dimensional Lie
algebras. More explicitly, if $Y$ is a smooth manifold acted upon
by a finite group $G$, then there exists a natural action on the
$n$-th Cartesian product $Y^n$ by the wreath product $\Gn$ (which
is the semidirect product of the symmetric group $S_n$ and the
product group $\G^n$). The orbifold Euler number for $Y^n/\Gn$ has
been explicitly calculated in {\em loc. cit.}. If in addition we
assume that $Y$ is a quasi-projective surface and $X$ is a
resolution of singularities of the orbifold $Y/G$, then the
following commutative diagram

$$\CD X^{[n]} @>>> X^n  /S_n \\ @VVV @VVV \\ Y^n  /G_n @<{\cong}<<
(Y/\G)^n /S_n \endCD $$
implies that the Hilbert scheme $X^{[n]}$ is a resolution of
singularities of the orbifold $X^n/ \Gn$. It has been shown
\cite{W} that if the ordinary Euler number of $X$ equals the
orbifold Euler number of $Y/\G$ then the ordinary Euler number of
$X^{[n]}$ equals the orbifold Euler number of $X^n/ \Gn$ for all
$n$. When $G$ is the trivial group and $X$ equals $Y$, one
recovers the case of symmetric products.

The purpose of the present paper is to point out that the wreath
product orbifolds also provide a large class of new higher
dimensional examples which verify the orbifold Hodge number
conjecture. More precisely, we show that if $Y$ is a
quasi-projective surface and $X$ is a resolution of singularities
of $Y/G$ such that the ordinary Hodge numbers of $X$ coincide with
the orbifold Hodge numbers of the orbifold $Y/G$, then the
orbifold Hodge numbers of the orbifold $Y^n/\Gn$ coincide with the
Hodge numbers of the Hilbert scheme $X^{[n]}$, which is a
resolution of singularities. Our proof rely on the analysis of
fixed-point set structures of the wreath product action on $Y^n$
(cf. \cite{W}). As in \cite{Zh} our calculation of the orbifold
Hodge numbers for $Y^n/\Gn$ actually works for any complex
$G$-manifold of even dimension.

In a very recent paper \cite{BDL}, Bryan, Donagi and Leung pointed
out a series of examples (besides the well-known symmetric
products) verifying the orbifold Hodge numbers conjecture. It
turns out that their examples correspond to our special case when
$Y$ is an abelian surface, $\G$ is $\Z_2$, and $X$ is the Kummer
$K3$ surface. They remarked that most examples in literature are
lower dimensional and no other higher dimensional examples known
to them. Their remarks are largely responsible for us to decide to
write up the results on the wreath product orbifolds which have
been known to us for some time. We do not know any other higher
dimensional examples which verify the orbifold Hodge number
conjecture.

To conclude we also make two explicit conjectures on elliptic
genera for wreath product orbifolds. These are motivated and in
turn generalize the work of Dijkgraaf, Moore, Verlinde and
Verlinde \cite{DMVV} on elliptic genera for symmetric products.

The layout of this paper is as follows. In
section~\ref{sec_prelim}, we recall the fixed-point set structures
of the wreath product action, and the definition of the orbifold
Hodge numbers, for both the compact and noncompact situations. In
section~\ref{sec_main}, we prove our main results,
Theorem~\ref{thm:main} on the orbifold Hodge numbers of wreath
product orbifolds and Theorem~\ref{samehodge} on the verification
of orbifold Hodge number conjecture. In section~\ref{sec_example},
we provide various examples illustrating our main results and in
addition formulate two conjectures on elliptic genera.

\section{preliminaries on the wreath product and orbifold Hodge
numbers} \label{sec_prelim}

In this section, we first review the definition of a wreath
product $\Gn$ associated to a finite group $\G$, and the
descriptions of conjugacy classes and centralizers for $\Gn$, cf.
e.g. \cite{Ker, M}. We also describe fixed-point sets for the
action of $\Gn$ on the $n$-th Cartesian product of a $\G$-manifold,
following \cite{W}. We then recall the definition of orbifold Hodge
numbers, cf. \cite{Zas, Bat-Dai}.
\subsection{The wreath product action on $Y^n$}
\label{sec:fixed}

Let $\G$ be a finite group and denote by $\G_*$ the set of
conjugacy classes of $\G$. Let $\G^n = \G \times \ldots \times \G$
be the direct product of $n$ copies of $\G$. Denote by $[g]$ the
conjugacy class of $g \in \G$. The symmetric group $S_n$ acts on
$\G^n$ by permuting the $n$ factors: $ s (g_1, \ldots, g_n) =
(g_{s^{-1}(1)} , \ldots, g_{s^{-1}(n)} )$. The {\em wreath
product} $\Gn = \G \wr S_n$ is defined to be the semidirect
product $\G^n \rtimes S_n$ of $\G^n$ and $S_n$, namely the
multiplication on $\Gn$ is given by $(g, s)(h, t) = (g. s(h),
st)$, where $g, h \in G^n, s, t \in S_n$. Note when $\G$ is the
trivial one-element group the wreath product $\Gn$ reduces to
$S_n$, and when $\G$ is $\Z_2$ the wreath product $\Gn$ is the
hyperoctahedral group, the Weyl group of type $C$.

Given $a = (g, s) \in \Gn$ where $g = (g_1, \ldots, g_n)$, we
write $s \in S_n$ as a product of disjoint cycles: if  $z= (i_1,
\ldots, i_r)$ is one of them, the {\em cycle-product} $g_{i_r}
g_{i_{r-1}} \ldots g_{i_1} $ of $a$ corresponding to the cycle $z$
is determined by $g$ and $z$ up to conjugacy. For each $c \in
\G_*$ and each integer $r \geq 1$, let $m_r (c)$ be the number of
$r$-cycles in $s$ whose cycle-product lies in $c$. Denote by $\rho
(c)$ the partition having $m_r (c)$ parts equal to $r$ ($r \geq
1$) and denote by $\rho = ( \rho (c) )_{c \in \G_*}$ the
corresponding partition-valued function on $\G_*$. Note that $||
\rho || : = \sum_{c \in \G_*} |\rho (c)| = \sum_{c \in \G_*, r
\geq 1} r m_r (c) = n$, where $| \rho(c)|$ is the size of the
partition $\rho (c)$. Thus we have defined a map from $\Gn$ to
${\mathcal P}_n (\G_*)$, the set of partition-valued function
$\rho =(\rho(c))_{c \in \G_*}$ on $\G_*$ such that $|| \rho ||
=n$. The function $\rho$ or the data $\{m_r(c) \}_{r,c}$ is called
the {\em type} of $a = (g, s) \in \Gn$. Denote ${\mathcal P}
(\G_*) = \sum_{n \geq 0} {\mathcal P}_n (\G_*)$. It is well known
(cf. e.g. \cite{Ker, M}) that two elements in $\Gn$ are conjugate
to each other if and only if they have the same type.

Let us describe the centralizer $Z_{\Gn} (a)$ of $a \in \Gn$, cf.
\cite{Ker, M, W}. First we consider the typical case that $a$ has
one $n$-cycle. As the centralizers of conjugate elements are
conjugate subgroups, we may assume that $a$ is of the form $a = (
(g, 1, \ldots, 1), \tau)$, where $ \tau = (1 2 \ldots n)$. Denote
by $Z_\G^{\Delta}(g)$, or $Z_\G^{\Delta_n}(g)$ when it is
necessary to specify $n$, the following diagonal subgroup of $G^n$
(and thus a subgroup of $\Gn$):
\begin{eqnarray*}
Z_G^{\Delta}(g) = \left\{ ( (h, \ldots, h), 1) \in G^n \mid h \in Z_G(g)
                 \right\}.
\end{eqnarray*}
The centralizer $Z_{\Gn} (a)$ of $a$ in $\Gn$ is equal to the
product $Z_G^{\Delta}(g) \cdot \langle a \rangle$, where $\langle
a \rangle$ is the cyclic subgroup of $\Gn$ generated by $a$.

Take a generic element $a  = ( g, s) \in \Gn$ of type $\rho = (
\rho (c) )_{c \in G_*}$, where $\rho (c) $ has $m_r (c)$
$r$-cycles ($r \geq 1$). We may assume (by taking a conjugation if
necessary) that the $m_r (c)$ $r$-cycles are of the form
$$
 g_{ur}(c) = ( (g, 1, \ldots,1), (i_{u1}, \ldots, i_{ur}) ),\quad
  1 \leq u \leq m_r (c), g \in c.
$$
Denote $ g_r (c) = ( (g, 1, \ldots,1), (12 \ldots r ) ).$
Throughout the paper, $\prod_{c,r}$ is understood as the product
$\prod_{c \in \G_*, r \geq 1}$. The centralizer $ Z_{\Gn} (a)$ of
$a \in \Gn$ is isomorphic to a direct product of the wreath
products
 \begin{eqnarray*}  \label{eq_centra}
  \prod_{c,r} \left(
                 Z_{G_r} ( g_r (c) )
     \wr S_{m_r (c)} \right).
 \end{eqnarray*}
Furthermore $Z_{G_r} ( g_r (c) )$ is isomorphic to $Z^{\Delta_r}_G
(g ) \cdot \langle g_{r} (c) \rangle$.

For a $\G$-space $Y$, we define an action of $\Gn$ on $Y^n$ as
follows. Given $ a = ( (g_1, \ldots, g_n), s)$, we let
\begin{eqnarray} \label{eq_action}
 a . (x_1, \ldots, x_n)
  = (g_1 x_{s^{-1} (1)}, \ldots, g_n x_{s^{-1} (n)})
\end{eqnarray}
where $x_1, \ldots, x_n \in Y$.

Next we recall the description of the fixed point set $( Y^n )^a$
for $a \in \Gn$, cf. \cite{W}. Let us first look at the typical
case $a = ( (g, 1, \ldots, 1), \tau) \in \Gn$. Note that the
centralizer group $Z_G(g)$ preserves the $g$-fixed point set
$X^g$. The fixed point set is
 \begin{eqnarray*}
  ( Y^n )^a = \left\{ (x, \ldots, x) \in Y^n\mid x= g x
              \right\}
 \end{eqnarray*}
which can be naturally identified with $Y^g$. The action of
$Z_{\Gn} (a)$ on $(Y^n)^a$ can be identified canonically with that
of $Z_\G (g)$ on $Y^g$ together with the trivial action of the
cyclic group $\langle a \rangle$. Thus $ (X^n)^a / Z_{\Gn} (a)$
can be identified with $ X^g / Z_\G (g).$

All $Z_\G (g)$ are conjugate and all $X^g$ are homeomorphic to
each other for different representatives $g$ in a fixed conjugacy
class $c \in G_*$. Also the orbit space $X^g /Z_G (g)$ can be
identified with each other by conjugation for different
representatives of $g$ in $c \in G_*$. We agree to denote $Z_\G
(g)$ (resp. $Y^g$, $Y^g /Z_\G (g)$) by $Z_\G (c)$ (resp. $Y^c$,
$Y^c /Z_\G (c)$) by abuse of notations. Similar remarks apply to
other situations below when the choice of representatives in a
conjugacy class is irrelevant.

For an element $a \in \Gn$ of type $\{m_r(c) \}$, the fixed-point
set $(Y^n )^a$ can be naturally identified with $\prod_{c,r} (Y^c
)^{m_r (c)} $. Furthermore the orbit space $(Y^n )^a /Z_{\Gn}(a)$
can be naturally identified with
\begin{eqnarray} \label{eq_fixedpoint}
 \prod_{c,r} S^{m_r (c)} \left( Y^{c} / Z_G (c) \right)
\end{eqnarray}
 where $S^{m}(X)$ denotes the $m$-th symmetric product $X^m/S_m$.
\subsection{Definition of orbifold Hodge numbers}

Let $Y$ be a compact complex manifold of complex dimension $d$
acted on by a finite group $G$ of automorphisms. For each
conjugacy class $c =[g] \in \G_*$, let $Y^g_1, \cdots, Y^g_{N_c}$
be the connected components of the fixed-point set $Y^g$. Zaslow
\cite{Zas} defined a shift number $F^g_{\alpha}$ associated to
each component $Y^g_{\alpha}$ as follows. On the tangent space to
each point in
 $Y^g_{\alpha}$,
$g$ acts as a diagonal matrix $\diag(e^{ 2\pi\sqrt{-1} \theta_1},
\cdots, e^{2\pi\sqrt{-1}\theta_d})$, where $0 \leq \theta_i <1$.
Then $$F^g_{\alpha} =\sum_{j=1}^d \theta_j.$$ In general,
$F^g_{\alpha}$ is just a rational number. However, there are many
occasions when it is an integer, e.g., when $g$ acts on the
tangent space by a matrix in $SL(n, \bC)$.

\begin{remark} \label{rem_integer}
In the case when $Y$ is a complex surface, the shift
$F_{\alpha}^g$ is an integer only if the component $Y^g_{\alpha}$
is either an isolated point or two dimensional. Indeed a finite
subgroup $\G$ of $GL(2, \C)$ acting on $\C^2$ has integer shifts
if and only if $\G$ lies in $SL(2, \C)$. That is, the shift
$F^g_{\alpha} = \theta_1 + \theta_2 $ is an integer if and only if
$\det g =e^{2\pi \sqrt{-1}(\theta_1 + \theta_2)} =1.$
\end{remark}

In the case all the shifts are integers, the {\em orbifold Hodge
numbers} of the orbifold $Y/\G$ are defined to be
\begin{eqnarray} \label{def_orbifold}
 h^{s,t}(Y, G) = \sum_{c \in G_*}
\sum_{\alpha_c =1}^{N_c} h^{s- F^c_{\alpha_c},
t-F^c_{\alpha_c}}(Y^c_{\alpha_c}/Z(c)).
\end{eqnarray}
The ordinary Dolbeault cohomology for an orbifold is given by (cf.
Satake \cite{Sat})

\begin{eqnarray} \label{def:OrbifolddeRham}
H^{*, *}(Y/G) \cong H^{*,*}(Y)^G.
\end{eqnarray}
Clearly the orbifold Hodge numbers can now be regarded as the
dimensions of the corresponding {\em orbifold cohomology groups}
(cf. \cite{Zas, Zh})

\begin{eqnarray} \label{str_cohomology}
 H^{*,*}(Y, G) = \sum_{c \in G_*} \bigoplus_{\alpha_c =1}^{N_c}
 H^{*,*}(Y^c_{\alpha_c}/Z_c)\{ F^c_{\alpha_c}\}.
\end{eqnarray}
Here and below we adopt the convention that if $V = \oplus_{s, t
\in \bZ} V^{s, t}$ is a bigraded vector space, then $V\{n\}$ is
the bigraded vector space with $(V\{n\})^{s, t} = V^{s-n, t-n}.$

It is convenient to form the generating function of bigraded
spaces

$$ H(Y,\G;x,y) =\sum_{s,t} H^{s,t}(Y,G) x^s y^t,$$
whose graded dimension is given by the orbifold Hodge polynomial
$$h(Y,\G;x,y) =\sum_{s,t} h^{s,t}(Y,G) x^s y^t.$$
Then we can rewrite
the definition of orbifold cohomology groups as

\begin{eqnarray} \label{eq_hodgepoly}
 H(Y,\G;x,y)
   &=& \sum_{c \in G_*} \bigoplus_{\alpha_c =1}^{N_c}
 H(Y^c_{\alpha_c}/Z_c; x,y)\{ F^c_{\alpha_c} \}  \\
   &=& \sum_{c \in G_*} \bigoplus_{\alpha_c =1}^{N_c}
 H(Y^c_{\alpha_c}/Z_c; x,y)(xy)^{ F^c_{\alpha_c}}.  \nonumber
\end{eqnarray}
For later use we
define the {\em orbifold virtual Hodge polynomial}

$$e(Y,\G;x,y) =\sum_{s,t} (-1)^{s+t}h^{s,t}(Y,G) x^s y^t.$$
We also define the usual virtual Hodge polynomial for the Hodge
numbers $h^{s,t}(Y)$ associated to smooth $Y$ by letting $e(Y;x,y)
=\sum_{s,t} (-1)^{s+t} h^{s,t}(Y) x^s y^t.$
\subsection{The definition of orbifold virtual Hodge numbers}

We now indicate how to extend the above definitions to the case of
smooth quasi-projective varieties by using Deligne's theory of
mixed Hodge structures \cite{Del}. Recall that a (pure) {\em Hodge
structure} of weight $m$ on a complex vector space $H$ with a real
structure is a direct sum decomposition: $$H = \bigoplus_{s + t =
m} H^{s, t},$$ such that $\overline{H}^{s, t} = H^{t, s}$ for all
pairs $(s, t)$. A {\em mixed Hodge structure} (MHS) on $H$
consists of two filtrations $$0 \subset \cdots \subset W_{m-1}
\subset W_m \subset W_{m+1} \subset \cdots \subset H,$$ the
`weight filtration', and $$H \supset  \cdots \supset F^{p-1}
\supset F^p \supset F^{p+1} \supset \cdots \supset 0,$$ the `Hodge
filtration', such that the filtration induced by the latter on
$Gr_m(W_*) = W_m/W_{m-1}$ defines a Hodge structure of weight $m$,
for each $m$.

Define
$$I^{s, t} = F^s\cap W_{s+t} \cap\left[ \overline{F^t} \cap W_{s+t}
 + \sum_{i\geq 2} \overline{F^{t-i+1}} \cap W_{s+t - i}\right];$$
Then $I^{s,t} \subset W_{s+t}$ maps isomorphically to the $(s,t)$ component
in $Gr_{s+t}(W_*)$.
One can show that
\begin{align*}
F^s(H) & = \bigoplus_{s' \geq s} \bigoplus_t I^{s', t}(H), &
W_m(H) & = \bigoplus_{s+t \leq m} I^{s, t}(H).
\end{align*}
It can be shown that $\{I^{s, t}\}$ is a splitting of $H$
characterized by the property that

$$I^{s, t} \cong \overline{I^{t, s}} \left( \mod \bigoplus_{s' <
s, t' < t} I^{s', t'} \right)$$ (cf. \cite{Del}). We will refer to
this splitting as the {\em canonical splitting}. Define
$$h^{s,t}(H) = \dim I^{s,t}(H).$$ Let $V = \oplus_{k \geq 0} V^k$
be a graded vector space, with $\dim V^k < \infty$ for all $k$.
Assume that each $V^k$ is endowed with a MHS. We will refer to
such a space as a {\em graded vector space with MHS}. The {\em
virtual Hodge numbers} and the {\em virtual Hodge polynomial} of
$V$ are defined by
\begin{eqnarray*}
&& e^{s, t}(V) = \sum_{k \geq 0} (-1)^k h^{s, t}(V^k), \\
 && e_{x,y}(V) = \sum_{s, t} e^{s, t}(V) x^sy^t.
\end{eqnarray*}
Alternatively, we have the splitting: $$V = \bigoplus_{k \geq 0}
\bigoplus_{s, t} I^{s, t}(V^k).$$ Consider the generating function
$$e_{x, y, z}(V) = \sum_{k \geq 0} \sum_{s, t} \dim I^{s, t}(V^k)
x^sy^tz^k.$$ Then $e_{x,y}(V) = e_{x, y,-1}(V)$. We will use the
following convention: for a graded vector space with MHS $V =
\oplus_{k \geq 0}V^k$ and a positive integer $n$, $V\{n\}$ is the
graded vector space with MHS such that for each $k$,
\begin{align*}
W_m(V^k \{n\}) & = W_{m - 2n}(V^k \{n\}) , & F^p(V^k \{n\}) & =
F^{p-n}(V^k \{n\}).
\end{align*}
It is straightforward to see that $e^{s, t}(V\{n\}) = e^{s-n,
t-n}(V)$, and so $$e_{x, y}(V\{n\}) = (xy)^n e_{x, y}(V).$$

Deligne \cite{Del} has shown that for an arbitrary complex
algebraic variety $Y$, the cohomology $H^k(Y)$ carries a MHS which
coincides with the classical pure Hodge structure in the case of
smooth projective varieties. Hence one can define the {\em virtual
Hodge number} of $Y$ $$e^{s, t}(Y) = e^{s, t}(H^*(Y))$$ and the
{\em virtual Hodge polynomial} of $Y$ $$e(Y; x, y) = e_{x,
y}(H^*(Y)).$$

Assume that $(Y, G)$ is a pair consisting of a smooth
quasi-projective variety $Y$ and a finite subgroup $G$ of
automorphisms of $Y$. Then by functorial property, there is an
induced action of $G$ on the MHS on $H^*(Y)$ by automorphisms. By
taking the invariant parts, we obtain a MHS on each $H^k(Y/G)$.
One can also achieve this by taking a smooth compactification
$\overline{Y}$ such that $D = \overline{Y} - Y$ is a divisor with
normal crossing singularities and such that the $G$-action extends
to $\overline{Y}$. Then the MHS on $H^*(Y/G)$ is obtained by using
$(\Omega^*_Y\langle D \rangle)^G$, the invariant part of the
complex of differential forms with logarithmic poles. Using the
above MHS on $H^*(Y/G)$, we can now define  $e^{p, q}(Y/G)$.
Similar to the closed case (cf. (\ref{def_orbifold})), we define
the {\em orbifold virtual Hodge number} as follows:
\begin{eqnarray*}
e^{s, t}(Y,G) = \sum_{c \in G_*} \sum_{\alpha_c = 1}^{N_c}
e^{s - F^c_{\alpha_c}, t - F^c_{\alpha_c}}(Y^c_{\alpha_c}/Z(c)).
\end{eqnarray*}
We also define the {\em orbifold virtual Hodge polynomial}:
\begin{eqnarray*}
e(Y, G; x, y) = \sum_{s, t} e^{s, t}(Y, G) x^sy^t.
\end{eqnarray*}
It is clear that $e(Y, G; x, y)$ is the virtual Hodge polynomial
of
$$H^*(Y, G) = \sum_{c \in G_*}
\bigoplus_{\alpha_c = 1}^N H^*(Y^c_{\alpha_c}/Z(c))\{F^c_{\alpha_c}\}$$
(cf (\ref{str_cohomology})),
where both sides are understood as graded vector spaces with MHS.

\begin{remark}
One can replace $H^*(Y)$ by the cohomology with compact support
$H^*_c(Y)$ in the above definitions.
\end{remark}

\section{The orbifold Hodge numbers of wreath product orbifolds}
\label{sec_main}

In this section, we calculate explicitly the ordinary and orbifold
Hodge numbers of wreath product orbifolds $Y^n/\Gn$ associated to
an even-dimensional orbifold $Y/G$.
\subsection{Two simple lemmas}
Let $V = \oplus_{s,t \in \Z_+} V^{s, t}$ be a bigraded complex
vector space, such that $\dim V^{s,t} < \infty$ for all $s, t$,
where $\Z_+$ is the set of non-negative integers. We introduce the
generating function

$$h_{x,y}(V) = \sum_{s, t \in \Z_+} (\dim V^{s, t}) x^s y^t.$$ For
example, when $V$ is the total Dolbeault cohomology group $H(Y)$,
then $h_{x,y}(V)$ is its associated Hodge polynomial $h(Y;x,y)$.
When $V$ is the total orbifold Dolbeault cohomology group
$H(Y,\G)$, then $h_{x,y}(V)$ is its associated orbifold Hodge
polynomial $h(Y,G;x,y)$. It is actually more convenient to work
with $e_{x, y}(V) = h_{-x, -y}(V)$. It is easy to see that
\begin{eqnarray*}
&& e_{x,y}(V_1
\oplus V_2) =
 e_{x,y}(V_1) + e_{x,y}(V_2), \\
&& e_{x,y}(V_1
\otimes V_2) =
 e_{x,y}(V_1) e_{x,y}(V_2).
\end{eqnarray*}
The graded symmetric algebra of $V$ is
 by definition $$S(V) = T(V)/I$$
where $T(V)$ is the tensor algebra
 of $V$,
$I$ is the ideal generated by elements of the form $$v
\otimes w - (-1)^{(s+t)(p+q)} w \otimes v, \;\; v \in V^{s,t}, \;
w \in V^{p, q}.$$ The bigrading on $V$ induces a bigrading
on $T(V)$ and also on
 $S(V)$,
and hence $e_{x,y}(S(V))$ makes sense. Note that for
bigraded vector spaces $V_1$ and $V_2$, we have $ S(V_1 \oplus
V_2) \cong S(V_1) \otimes S(V_2)$.
Consequently,

\begin{eqnarray} \label{eq_multiplicative}
e_{x,y} ( S(V_1 \oplus V_2)) = e_{x,y}(S(V_1))
e_{x,y}(S(V_2)).
\end{eqnarray}
By introducing a formal variable $q$ to count the degree of
symmetric power, we can write formally $ S(qV) = \sum_{n \geq 0}
S^n(V) q^n$. By breaking $V$ into one-dimensional subspaces, one
can easily prove the following.

\begin{lemma} \label{lm:Hodge}
For any bigraded vector space $V=\oplus_{s,t \geq 0} V^{s,t}$
with $\dim V^{s, t} < \infty$ for all pairs $(s, t)$,
we
 have
$$\sum_{n \geq 0} e_{x,y}(S^n(V)) q^n = \prod_{s, t}
\frac{1} { (1 -
 x^sy^tq)^{e^{s, t}(V)}},$$
where $e^{s, t}(V) = (-1)^{s+t}\dim V^{s, t}$.
\end{lemma}

For a formal power series $\sum_{r > 0}V_r q^r$, where each $V_r$
is a bigraded vector space of weight $r$ such that $\dim V_r^{s,
t} < \infty$, define

$$S(\sum_{r> 0}V_r q^r) = \sum_{m \geq 0} \sum_{\sum_{j=1}^m j m_j
= m} q^m \bigotimes_{j=1}^m S^{m_j}(V_j).$$

Formally we have

$$S(\sum_{r > 0}V_r q^r) = \bigotimes_{r> 0} S(V_r q^r)$$ and

$$e_{x,y}(\sum_{r > 0}V_r q^r) = \sum_{r> 0} e_{x,y}(V_r) q^r.$$
Then the next lemma follows from Lemma \ref{lm:Hodge}.

\begin{lemma} \label{lm:Hodge2}
For a sequence $\{V_n\}$ of bigraded vector spaces,
we have the following formula:
$$e_{x,y}\left( S(\sum_{n > 0}
 V_n q^n) \right)
= \prod_{n > 0} \prod_{s, t}\frac{1} {(1 -
 x^sy^tq^n)^{e^{s, t}(V_n)}}.$$
\end{lemma}

\begin{remark} \rm  \label{rem_mixhodge}
Using the canonical splitting, it is fairly straightforward to
generalize Lemma \ref{lm:Hodge} and Lemma \ref{lm:Hodge2} to the
case of vector spaces with MHS.
\end{remark}

\subsection{The main theorems}

Since $\G^n$ is a normal subgroup of the wreath product $\Gn =\G^n
\rtimes S_n$, it is easy to see by (\ref{def:OrbifolddeRham}) that

\[
H(Y^n/\Gn; x,y) \cong H(Y^n;x,y)^{\G^n \rtimes S_n} \cong S^n(H(Y;
x,y)^G) \cong S^n(H(Y/G; x,y)).
\]
When $Y$ is a compact complex manifold, this is an isomorphism of
bigraded vector spaces; when $Y$ is a quasi-projective smooth
variety over $\bC$, this is an isomorphism of graded vector spaces
with MHS. As a consequence of Lemma~\ref{lm:Hodge} and
Remark~\ref{rem_mixhodge}, we obtain the following proposition.

\begin{proposition}
If $Y$ is a compact complex manifold or a quasi-projective smooth
variety, and $G$ is a finite subgroup of automorphisms on $Y$,
then we have the following formula:

\begin{eqnarray*}
\sum_{n \geq 0} e(Y^n/\Gn; x,y) q^n = \prod_{s, t}
\frac{1} {(1 - x^sy^t
 q)^{e^{s,t}(Y/G)}}.
\end{eqnarray*}
\end{proposition}

The first main result of this paper is the following theorem.

\begin{theorem}\label{thm:main}
Given a compact complex manifold or a smooth quasi-projective
variety $Y$ of even complex dimension $d$,
acted on by a finite group $G$ with integer shifts,
% for the wreath product orbifold $Y^n/\Gn$ is given by
we have the
 following formula for the orbifold Hodge
 numbers:

\begin{eqnarray} \label{eqn:generating}
\sum_{n=1}^{\infty} e(Y^n, \Gn; x,y) q^n =
\prod_{r=1}^{\infty}\prod_{s,t}\frac{1} { (1 - x^sy^t
q^r(xy)^{(r-1)d/2})^{e^{s, t}(Y, G)}} .
\end{eqnarray}
%Here we have assumed that all the shifts for $(Y, G)$ are integers.
\end{theorem}

\begin{proof}
We first compute the shifts $F^c$ for the orbifold $Y^n/\Gn$
associated to a conjugacy class $c$ in $\Gn$. Consider the typical
class containing $$g \wr \sigma_n =((g, 1, \cdots, 1), (12\cdots
n))$$ where $ g \in \G$. Recall from the previous section that a
fixed point in $Y^n$ by the action of $g \wr \sigma_n$ is of the
form $(x, \ldots, x)$ where $x \in Y^g$. Since the calculation can
be done locally, we will assume that we take local coordinates
$(z_1, \cdots, z_d)$ near a point $x \in Y^g$ such that the action
is given by

$$g (z_1, \cdots, z_d) = (e^{2\pi\sqrt{-1}\theta_1}z_1, \cdots,
e^{2\pi\sqrt{-1}\theta_r}z_r, z_{r+1}, \cdots, z_d).$$
Equivalently, $g$ is locally given by the diagonal matrix
$\diag(e^{2\pi\sqrt{-1}\theta_1}, \cdots,
e^{2\pi\sqrt{-1}\theta_d})$ where $\theta_{r+1} = \cdots =
\theta_d = 0$. Then on $Y^n$ near $(x, \cdots, x)$, $g \wr
\sigma_n$ is given by a block diagonal matrix with blocks of the
form

$$\left(
\begin{array}{ccccc} 0 & 1 & 0 & \cdots & 0 \\ 0 & 0 & 1 & \cdots
& 0 \\ \cdots & \cdots & \cdots & \cdots & 0 \\ 0 & 0 & 0 & \cdots
& 1 \\ e^{2\pi\sqrt{-1}\theta_j} & 0 & \cdots & \cdots & 0
\end{array} \right)$$

The characteristic polynomial of this matrix is $\lambda^n -
e^{\sqrt{-1}\theta_j}$, hence it has eigenvalues $$\lambda_{jk} =
e^{2\pi \sqrt{-1}(\theta_j+k)/n}, \;\; k =0, \cdots, n-1.$$ Notice
that $\lambda_{jk} = 1$ if and only $\theta_j = k = 0$. So the
shift for the component of $(Y^n)^{g \wr \sigma_n}$ containing
$(x, \ldots, x)$ is given by

\begin{eqnarray*}
F^{g\wr\sigma_n}(x, \ldots, x)
 & =& \sum_{j=1}^r \sum_{k=0}^{n-1} \frac{\theta_j+k}{n} + (d-r)
\sum_{k=1}^{n-1}\frac{k}{n} \\
 & =&
\sum_{j=1}^r \theta_j + (n-1)d/2 = F^c_{\alpha_c} + (n-1)d/2.
\end{eqnarray*}
Here we have assumed that $x \in Y^g$ lies in the component
$Y^g_{\alpha_c}$ $({\alpha_c} =1, \ldots, N_c)$, and
$F^c_{\alpha_c}$ is the shift for the component
$Y^c_{\alpha_c}/Z_{\G}(c)$.

For a general conjugacy class containing an element $a$ of type
$$\rho =\{m_r(c)\}_{r \geq 1, c \in \G_*},$$ where $\sum_{r, c}r
m_r(c) =n$, the description of the fixed-point set $(Y^n)^a$ given
in (\ref{eq_fixedpoint}) implies that the components for $(Y^n)^a$
can be listed as
 $$
 (Y^n)^a_{\{m_{r,c}(\alpha_c)\}} =
 \prod_{r,c} \prod_{\alpha_c =1}^{N_c} S^{m_{r,c}(\alpha_c)}
 (Y^c_{\alpha_c} /Z_\G(c)),
 $$
where $(m_{r,c}(1), \ldots, m_{r,c}(N_c))$ satisfies
$\sum_{\alpha_c =1}^{N_c} m_{r,c}(\alpha_c) = m_{r}(c)$. Then the
shift for the component $(Y^n)^a_{\{m_{r,c}(\alpha_c)\}}$ is given
by
\begin{eqnarray} \label{eq_shiftnumber}
F_{\{m_{r,c}(\alpha_c)\}}
 = \sum_{r , c} \sum_{\alpha_c
=1}^{N_c} m_{r,c}(\alpha_c) \left( F^c_{\alpha_c} + (r-1)d/2
\right).
\end{eqnarray}

By using (\ref{eq_hodgepoly}), (\ref{eq_fixedpoint}),
(\ref{eq_shiftnumber}) and (\ref{eq_multiplicative}) we have

\begin{eqnarray*}
&& \sum_{n \geq 0} H(Y^n, \Gn; x,y)q^n \\
 & = &
 \sum_{n =0}^{\infty}
 \bigoplus_{\{m_r(c) \} \in {\mathcal P}_n (\G_*)}
 \bigotimes_{r,c} \bigotimes_{\alpha_c =1}^{N_c}
 H(S^{m_{r,c}(\alpha_c)}
 (Y^c_{\alpha_c} /Z_\G(c));x,y)\{F_{\{m_{r,c}(\alpha_c)\}} \} q^n \\
& = &
 \sum_{n =0}^{\infty}
 \bigoplus_{\{m_r(c) \} \in {\mathcal P}_n (\G_*)} \\
 & & \qquad \bigotimes_{r,c} \bigotimes_{\alpha_c =1}^{N_c}
 S^{m_{r,c}(\alpha_c)}
 \left(H (Y^c_{\alpha_c} /Z_\G(c);x,y)
 \{ F^c_{\alpha_c} + (r-1)d/2 \}
 \right) q^n \\
& = &
 \sum_{\{m_r(c) \}}
 \bigotimes_{r,c}  S^{m_{r,c}}
 \left( \bigoplus_{\alpha_c =1}^{N_c} H (Y^c_{\alpha_c} /Z_\G(c);x,y)
  \{ F^c_{\alpha_c} + (r-1)d/2 \} q^r\right)  \\
& = &
 \sum_{\{m_r \}}
 \bigotimes_{r}  S^{m_r}
 \left(\bigoplus_c \bigoplus_{\alpha_c =1}^{N_c}
 H (Y^c_{\alpha_c} /Z_\G(c);x,y)
  \{ F^c_{\alpha_c} + (r-1)d/2 \} q^r \right)  \\
& & \qquad\qquad\qquad\qquad\qquad\qquad\qquad \mbox{where we let
} m_r = \sum_c m_r(c) \\
 & = &
 \sum_{\{m_r \}} \bigotimes_{r \geq 1}  S^{m_r}
 \left( H (Y, \G;x,y)
\{ (r-1)d/2 \}  q^r\right)\\
 & = &
 S\left(\sum_{r \geq 1}  H(Y, \G;x,y)
\{ (r-1)d/2 \} q^r \right).
\end{eqnarray*}

Namely we have proved that

\begin{eqnarray*}
\sum_{n \geq 0} H(Y^n, \Gn; x,y)q^n = S \left( \sum_{r>0} H(Y,
G;x,y) (xy)^{ (r-1)d/2 } q^r \right),
\end{eqnarray*}
which implies immediately the theorem by means of
Lemma~\ref{lm:Hodge2}.
\end{proof}

\begin{remark} \rm
When $\G$ is trivial and $Y$ is an algebraic surface,
Theorem~\ref{thm:main} recovers the orbifold Hodge numbers for the
symmetric product $Y^n/S_n$ which was calculated in \cite{G2, Zh}.
On the other hand, if we set $x=y =1$ we recover the orbifold
Euler numbers for $Y^n/\Gn$ which was first computed in \cite{W}
for any topological space $Y$.
\end{remark}

\begin{remark} \rm
In the above we have constrained ourselves
to the case that the shift numbers are integers.
Physicists are also interested in the case of fractional shift numbers
(see e.g. Zaslow \cite{Zas}).
It is straightforward to generalize our result.
\end{remark}
\subsection{Some consequences}
We assume that $Y$ is a quasi-projective surface acted upon by a
finite group $\G$, and that $X$ is a resolution of singularities
of the orbifold $Y/\G$. We denote by $X^{[n]}$ the Hilbert scheme
of $n$ points on $X$. It is well known (cf. \cite{Fog, G1}) that
the Hilbert-Chow morphism $X^{[n]} \rightarrow X^n/S_n$ is a
resolution of singularities. Indeed it is crepant. We have the
following commutative diagram \cite{W}

$$\CD X^{[n]} @>>> X^n  /S_n \\ @VVV @VVV \\ Y^n  /G_n @<{\cong}<<
(Y/\G)^n /S_n \endCD $$
which implies that the Hilbert scheme $X^{[n]}$ is a resolution of
singularities of the orbifold $X^n/ \Gn$.

As calculated in \cite{Got-Soe} and \cite{Che}, the Hodge numbers
for the Hilbert scheme $X^{[n]}$ are given by the following
formula:

\begin{eqnarray*}
\sum_{n=1}^{\infty} e(X^{[n]}; x,y) q^n =
\prod_{r=1}^{\infty}\prod_{s, t}
\frac{1} {(1 - x^sy^t
 q^r(xy)^{r-1})^{e^{s, t}(X)}} .
\end{eqnarray*}
Here and below we use the cohomology with compact supports.

By comparing with Theorem~\ref{thm:main} we obtain the following
theorem which provides us a large class of higher dimensional
examples which verify the orbifold Hodge number conjecture (cf.
\cite{Vafa, Zas, Bat-Dai, Bat}). The assumption of the theorem is
necessary by Remark~\ref{rem_integer}.

\begin{theorem} \label{samehodge}
Let $Y$ be a smooth quasi-projective surface which admits a
$G$-action with only isolated fixed points. Assume that $\pi: X
\to Y/G$ is a resolution such that $e(X;x,y) = e(Y, G;x,y)$. Let
$X^{[n]}$ be the Hilbert scheme of $n$ points of $X$. Then for all
$r,s$ we have $$h^{r,s}(X^{[n]}) =h^{r, s}(Y^n, \Gn).$$
\end{theorem}

\begin{remark} \rm
When $\G$ is trivial and $X$ equals $Y$, we recover the theorem of
\cite{G2, Zh}. We will see later many interesting examples arise
when $\G$ is not trivial.
\end{remark}

More generally if $Y$ has dimension greater than $2$, there is no
such a favorable resolution as Hilbert scheme for $Y^n/\Gn$.
Nevertheless we have the following interesting corollary of
Theorem~\ref{thm:main}. Here we assume that the shifts are
integers for the orbifold $Y/\G$ so that its orbifold Hodge
numbers are well defined.

\begin{corollary} \label{cor1}
Let $Y$ be a smooth variety of even dimension
acted on by a finite group $G$ of automorphisms, and $\pi: X \to
Y/G$ is a resolution such that $h^{s, t}(X) = h^{s, t}(Y, G)$ for
all $s, t$, then for all $s, t$ we have $$h^{s, t}(X^n, S_n) =
h^{s, t}(Y^n, G_n).$$
\end{corollary}

\section{Examples and applications} \label{sec_example}

In this section we provide various concrete examples which satisfy
the assumptions of Theorem \ref{thm:main} and
Theorem~\ref{samehodge}. We also give explicit conjectures on the
elliptic genera for the wreath product orbifolds.

\subsection{Various examples}

\begin{example}
When $\G$ is trivial and $X$ equals $Y$, this gives us the example
of symmetric products \cite{G1, Zh}.
\end{example}

\begin{example}
$Y$ is $\C^2$, $\G$ is a finite subgroup of $SL_2(\C)$, and $X$ is
the minimal resolution of $\C^2/\G$. The
exceptional fiber consists of $|\G_*| -1$ irreducible components
which are $(-2)$-curves (cf. e.g. \cite{HH}). We have
\begin{eqnarray*}
h^{s,t}(X) = \left\{ \begin{array}{ll} 1, & s =t =0, \\ |\G_*| -1,
& s =t = 1, \\  0, & \text{otherwise}. \end{array} \right.
\end{eqnarray*}
On the other hand, for any non-trivial conjugacy class $c \in
\G_*$, the corresponding shift is $1$ and thus makes a
contribution to $h^{1,1}(\C^2, \G)$ which results that
$h^{1,1}(\C^2, \G) =|\G_*| -1.$ The other $h^{s,t}(\C^2, \G)$ can
be also seen to coincide with $h^{s,t}(X)$.

This example has played a key role in the connections between the
wreath product orbifolds and the vertex representations of affine
and toroidal Lie algebras \cite{W, FJW, W2}.
\end{example}

\begin{example} (Bryan-Donagi-Leung \cite{BDL})
Let $Y$ be an abelian surface (two dimensional torus). The
$\bZ_2$-action induced by the involution $\tau: x \to - x$ has
$16$ fixed points, at each of which the shift $F^{\tau}$ is $ 1$.
So the twisted sectors contribute an extra $16$ to $h^{1, 1}$.
Write $Y = \bC^2/L$ for some lattice $L$, and let $(z_1, z_2)$ be
the linear coordinates on $\bC^2$. Then $H^{*, *}(Y)$ is generated
by $d z^1, d\bar{z}^1, dz^2, d\bar{z}^2$. The action of $\tau$
just takes $dz^j$ to $-dz^j$, etc. Hence it is clear that

\begin{eqnarray*}
H^{*, *}(Y)^{\bZ_2} & \cong & \bC \oplus \bC dz^1 \wedge dz^2
\oplus (\oplus_{j, k =1}^2 \bC dz^j \wedge  d\bar{z}^k) \\
&& \oplus \bC d\bar{z}^1 \wedge d\bar{z}^2 \oplus \bC dz^1 \wedge
dz^2 \wedge d\bar{z}^1 \wedge d\bar{z}^2.
\end{eqnarray*}

Therefore,
\begin{eqnarray*}
h^{s,t}(Y, \bZ_2) = \left\{ \begin{array}{ll} 1, & s =t =0, \\ 20,
& s =t = 1, \\ 1, & s= 2, t = 0 \; \text{or} \; s=0, t =2, \\ 1, &
s =t= 2, \\ 0, & \text{otherwise}. \end{array} \right.
\end{eqnarray*}

The minimal resolution $X \to Y/\pm 1$ is a crepant resolution,
where $X$ is a K3 surface. This is the famous Kummer construction.
By the well known Hodge numbers of a K3 surface, we have $h^{s,
t}(X) = h^{s, t} (Y, \bZ_2)$ for all $s,t$.
\end{example}

\begin{example}
Let $\bZ_3$ act on $\bCP_2$ by

$$\alpha \cdot [z_0:z_1:z_2] = [\alpha z_0:\alpha^{-1}z_1:z_2],$$
where $\alpha$ is a generator of $\bZ_3$ and identified with a
cubic root of unity on the right-hand side. This action has three
fixed points: $p_0 = [1:0:0]$, $p_1=[0:1:0]$, and $p_2 = [0:0:1]$.
At these point, the weights of the action are $(1, 2)$, $(2, 1)$,
and $(1, 2)$ respectively. It is then straightforward to see that
for $g \neq 1$ we have

$$F^g = \frac{1}{3} + \frac{2}{3} = 1.$$
Therefore,
\begin{eqnarray*}
H^{*, *} (\bCP_2, \bZ_3) = H^{*, *}(\bCP_2)^{\bZ_3} \bigoplus
\bigoplus_{j=0}^2 H^{*, *}(p_j)^{\bZ_3}\{1\} \bigoplus
\bigoplus_{j=0}^2 H^{*, *}(p_j)^{\bZ_3}\{1\},
\end{eqnarray*}
and hence
\begin{eqnarray*}
h^{s,t}(\bCP_2, \bZ_3) = \left\{ \begin{array}{ll} 1, & s =t = 0,
\\ 7, & s =t = 1, \\ 1, & s =t = 2, \\ 0, & \text{otherwise}.
\end{array} \right.
\end{eqnarray*}
The minimal resolution $X = \widehat{\bCP_2/\bZ_3}$ is obtained by
replacing each singular point by a string of two $(-2)$-curves, each of
which contributes $1$ to $h^{1, 1}$, hence $h^{1, 1}$ of
$\widehat{\bCP_2/\bZ_3}$ is $7$. This resolution is a crepant
resolution.
\end{example}

\begin{example}
 \label{exm:A_n}
Let $n > 2$ be an odd number,
Consider the
action of $\bZ_n$ on $\bCP_3$ given by

$$\alpha \cdot [z_0:z_1:z_2:z_3]
= [z_0:z_1:\alpha
 z_2:\alpha^{-1}z_3],$$
where $\alpha$ is a generator of $\bZ_n$. It has a fixed line
$\{[z_0:z_1:0:0]\}$ and two isolated fixed points $[0:0:1:0]$, and
$[0:0:0:1]$. Let $Y_{m, n}$ be the Fermat surface defined by
$$z_0^{mn} + z_1^{mn} + z_2^{mn} + z_3^{mn} = 0$$ in $\bCP_3$. The
above action preserves $Y_{m, n}$, with $mn$ isolated fixed
points:

$$[1:e^{(2k +1) \pi \sqrt{-1}/(mn)}:0:0], \;\;\; k = 0, \dots,
mn-1.$$ Note the action is semi-free, i.e. the stabilizers are
either trivial or the whole group $\bZ_n$. Near each
 of the fixed points, say $[1: e^{\pi \sqrt{-1}/(mn)}:0:0]$,
$Y_{m,n}$ is given by the equation $$1 + u_1^n + u_2^n + u_3^n =
0,$$ where $u_j = z_j/z_0$. We can use $(u_2, u_3)$ as local
coordinates, then $\bZ_n$ acts with weight $(1, -1)$, i.e.,
$\bZ_n$ acts locally by matrices in $SL(2, \bC)$. Therefore,
$Y_{m, n}/Z_n$ admits a crepant resolution obtained by replacing
each isolated singular point with a string of $n-1$ copies of
$(-2)$-curves.
\end{example}

\begin{example}
 \label{exm:D_n}
Denote now by $\beta$ a generator of $\bZ_4$.
Consider the $\bZ_4$-action on $\bCP_3$ given by
$$\beta \cdot [z_0:z_1:z_2:z_3]
= [z_0:z_1:\sqrt{-1}
 z_3:\sqrt{-1}z_2].$$
Combined with the $\bZ_n$-action in Example \ref{exm:A_n}, we get
an action of the binary dihedral group $D_n^*$ on $\bCP_3$ which
preserves $Y_{4m, n}$. By the same method as in Example
\ref{exm:A_n} one can find the fixed points and sees that $Y_{4m,
n}/D_n^*$ admits a crepant resolution.
\end{example}

\begin{example}

The method of Example \ref{exm:A_n} and Example \ref{exm:D_n} can
be generalized to other finite subgroups of $SL(2, \bC)$. Given
such a group $G$, let it act on $\bC^4$ on the last two factors.
This action induces an action on $\bCP_3$. Now consider a smooth
hypersurface $Y$ defined by an equation of the form $$f(z_0, z_1)
+ g(z_2, z_3) = 0,$$ where $f$ and $g$ are two homogeneous
polynomials of the same degrees, and $g$ is an invariant
polynomial for $G$. Using the explicit description of the
$G$-action on $\bC^2$ and the invariant polynomials (see e.g.
Klein \cite{Kle}), one can find many examples which admits crepant
resolutions. One should be able to find more examples by
considering complete intersections in (weighted) projective
spaces.
\end{example}

\begin{example}
More complicated examples can be found in two papers by Barlow
\cite{Bar1, Bar2} , e.g. the quotient of a Hilbert modular surface
by $\bZ_2$ or $D_{10}$, or the quotient of a complete intersection
of $4$ quadrics in $\bCP_6$ by a group of order $16$, or the
quotient of a Godeaux-Reid surface by an involution.
\end{example}

\subsection{Conjectures on elliptic genera of wreath product orbifolds}

Let $Y$ be a compact K\"ahler manifold of complex dimension $d$,
denote by $TY$ (resp. $T^*Y$) its holomorphic tangent (resp.
cotangent) bundle. Consider the formal power series of vector
bundles:

$$E_{q, y}(Y) = y^{-\frac{d}{2}} \bigotimes_{n \geq 1} \left(
\Lambda_{-yq^{n-1}}(T^*Y) \otimes \Lambda_{-y^{-1}q^n} (TY)
\otimes S_{q^n}(T^*Y) \otimes S_{q^n}(TY)\right).$$
If we write

$$E_{q, y}(Y) = \sum_{m \geq 0,l} q^my^l E_{m, l}(Y),$$ we easily
see that each $E_{m, l}$ is a holomorphic bundle of finite rank,
hence one can consider its Riemann-Roch number $$c(m, l) =
\chi(E_{m, l}(Y)) = \sum_{k \geq 0} (-1)^k \dim H^k(Y,E_{m,
l}(Y)).$$ The generating function $$\chi(Y;q, y) = \sum_{m\geq 0,
l} q^my^l\chi(E_{m, l}(Y)) = \chi(E_{q, y}(Y))$$ is called the
{\em elliptic genus} of $Y$ (cf. \cite{Hir, Lan}). In the very
important special case when $q = 0$, one recovers the Hirzebruch
genus:
\begin{eqnarray*}
&& E_{0, y}(Y) = y^{-\frac{d}{2}}\Lambda_{-y}(T^*Y), \\
&& \chi(Y; 0, y) = y^{-\frac{d}{2}} \chi_{-y}(Y)
= y^{-\frac{d}{2}} \sum_{s, t \geq 0} (-1)^t(-y)^sh^{s, t}(Y).
\end{eqnarray*}

We do not know of a good mathematical formulation of elliptic
genera for orbifolds. However physicists have interpreted elliptic
genera as partition functions of supersymmetric sigma models,
which makes sense also for orbifolds (cf. \cite{Lan, DMVV} and
references therein). Based on physical arguments and the
description of fixed-point sets for the symmetric group action on
$Y^n$, Dijkgraaf {\em et al} \cite{DMVV} derived a formula for the
elliptic genera of the symmetric products $S^n(Y)$ in terms of
that of $Y$. In the case of a K3 surface or an abelian surface,
they also conjectured that the same formula should compute the
elliptic genera of the Hilbert schemes. Their method, if can be
made mathematically rigorous, should also provide the proof of the
following conjectures with suitable modifications.
%We make two conjectures on the elliptic
%genera for wreath product orbifolds which are parallel to our
%Theorem~\ref{thm:main} and Theorem~\ref{samehodge} (also see
%\cite{W} for the corresponding results on orbifold Euler numbers).
(In the following we denote by $\chi(Y, G; q, y)$ the elliptic genera
of an orbifold $Y/G$.
)

\begin{conjecture} \label{conj_wreath}
Let $Y$ be a K\"ahler $\G$-manifold. If we write the elliptic
genus for $Y/\G$ as $\chi (Y,G;q,y) = \sum_{m \geq 0, l} c(m,l)
q^m y^l$, then the elliptic genus for the wreath product orbifold
$Y^n/\Gn$ is given by the following formula:

$$\sum_{N=0}^{\infty} p^N \chi(Y^N, \G_N;q,y) = \prod_{n>0, m\geq
0,l}\frac1{(1 -p^nq^my^l)^{c(nm,l)}}. $$
\end{conjecture}

\begin{conjecture}
Let $Y$ be a K\"ahler $\G$-surface. We assume that $X$ is a
resolution of singularities of $Y/\G$ such that $\chi (Y,\G;q,y)
=\chi (X;q,y).$ Then $\chi (Y^n,\Gn;q,y) =\chi (X^{[n]};q,y)$ for
all $n$.
\end{conjecture}

When $G$ is trivial, one recovers the symmetric product situation
as in \cite{DMVV}. In this case, the $q = 0$ version of
Conjecture~\ref{conj_wreath} has been verified in \cite{Zh} as a
corollary of the calculation of orbifold Hodge numbers. Similarly,
our results in Section~\ref{sec_main} can be viewed as supporting
evidence for the above conjectures in the general setup of wreath
product orbifolds.

{\em Note added.} In a recent remarkable paper \cite{BL}, Borisov
and Libgober have introduced the mathematically rigorous notion of
orbifold elliptic genera among other things, and verified our
Conjecture~\ref{conj_wreath}.

\end{document}